\documentclass[12pt]{amsart}
\usepackage{amsfonts, amssymb, amsmath, amsthm,
latexsym}
\usepackage{a4wide}
\usepackage{amsrefs}
\usepackage{hyperref}

\begin{document}

\theoremstyle{plain}
\newtheorem*{main theorem}{Main Theorem}
\newtheorem*{theorem*}{Theorem}
\newtheorem*{emptytheorem*}{}
\newtheorem{theorem}{Theorem}
\newtheorem{lemma}{Lemma}
\newtheorem{sublemma}{Sublemma}[lemma]
\newtheorem{proposition}{Proposition}[section]
\newtheorem*{proposition*}{Proposition}
\newtheorem{corollary}{Corollary}[section]
\newtheorem{claim}{Claim}[sublemma]

\theoremstyle{definition}
\newtheorem{definition}{Definition}
\newtheorem*{definition*}{Definition}
\theoremstyle{remark}
\newtheorem{remark}{Remark}
\newtheorem*{remark*}{Remark}

\def \a{\alpha}
\newcommand{\g}{\gamma}
\def \l{\lambda}
\def \L{\Lambda}
\def \G{\Gamma}
\def \D{\Delta}
\def \d{\delta}
\def \e{\epsilon}
\def \g{\gamma}
\def \t{\theta}
\def \w{\omega}
\def \Om{\Omega}
\newcommand{\p}{\Phi}
\def \R{\mathbb{R}}
\def \Z{\mathbb{Z}}
\def \Boxx{\hfill$\blacksquare$}
\def \C{\EuScript{C}}
\def \E{\EuScript{E}}
\def \F{\EuScript{F}}
\def \I{\EuScript{I}}
\def \O{\EuScript{O}}
\def \A{\mathcal{A}}
\def \B{\mathcal{B}}
\def \K{\mathcal{K}}
\def \M{\mathcal{M}}
\def \N{\mathcal{N}}
\def \Q{\mathcal{Q}}
\def \T{\mathcal{T}}
\def \PP{\mathcal{P}}
\def \dd{\mathfrak{D}}

\newcommand{\mcal}[2]{\ensuremath{\mathcal #1^{(#2)}}}
\newcommand{\mcalp}[2]{\ensuremath{\mathcal #1^{(#2)}}(p)}

\newcommand{\pfi}[2]{\ensuremath{\partial_{#1}\Phi_{#2}}}
\newcommand{\pd}[2]{\ensuremath{\partial_{#1}^{#2}}}
\newcommand{\pdtil}[2]{\ensuremath{\tilde\partial_{#1}^{#2}}}

\newcommand{\coni}[1]{\ensuremath{e^{(#1)} }}
\newcommand{\exi}[1]{\ensuremath{f^{(#1)} }}
\newcommand{\con}[2]{\ensuremath{e^{(#1)}(#2) }}
\newcommand{\ex}[2]{\ensuremath{f^{(#1)}(#2) }}
\newcommand{\conim}[3]{\ensuremath{e^{(#1)}_{#2}(#3) }}
\newcommand{\exim}[3]{\ensuremath{f^{(#1)}_{#2}(#3) }}
\newcommand{\conimi}[2]{\ensuremath{e^{(#1)}_{#2} }}
\newcommand{\eximi}[2]{\ensuremath{f^{(#1)}_{#2} }}

\author{Mark Holland}
\address{Mathematics Department, University of Surrey, Guildford GU2
SXH}
\email{mark.holland@eim.surrey.ac.uk}
\author{Stefano Luzzatto}
\address{Mathematics Department, Imperial College, London SW7 2AZ}
\email{Stefano.Luzzatto@Imperial.ac.uk}
\urladdr{http://www.ma.ic.ac.uk/\textasciitilde luzzatto}
\title[Stable manifolds under
very weak hyperbolicity conditions]{Stable manifolds under
very weak \\ hyperbolicity conditions}

\date{\today}

\begin{abstract}
We present an argument for proving the existence of local stable
and unstable manifolds in a general abstract setting and under 
very weak hyperbolicity conditions.    
\end{abstract}

\maketitle

\section{Introduction and Results}

\subsection{Stable sets and stable manifolds}

One of the most fundamental concepts in the modern geometric theory 
of dynamical systems is that of the \emph{stable set} associated to a 
point: given a map \( \varphi: M \to M \) on a metric space
\( M \), and a point \( 
z\in M \) we define the (global) \emph{ stable set} of \( z \) as
    \[ 
    W^{s}(z) = \{x\in M : 
    d(\varphi^{k} (x), 
    \varphi^{k}(z)) \to 0 \text{ as } k\to \infty\}.
    \]
If \( \varphi \) is invertible,  (global) \emph{ unstable set} can be defined 
in the same way by taking \( k\to -\infty \). The situation is completely 
analogous and so we will concentrate here on stable sets. 

This definition gives an equivalence relation on \( M \) which 
defines a partition into  sets 
which are invariant under the action of \( \varphi \) and which are 
formed of orbits which have the same asymptotic behaviour. An 
understanding of the geometry of the stable and unstable sets of 
different points, of how they depend on the base point \( z \), and of
how they intersect, forms the  
core of many powerful arguments related to all kinds of 
properties of dynamical systems, from ergodicity to structural 
stability to estimates on decay of correlations. 

In general \( W^{s}(z) \) can be 
extremely complicated, both in its intrinsic geometry and/or
in the way it is embedded in \( M \). A first step towards 
understanding this complexity is to focus on the 
\emph{local stable set}
    \[ 
    W^{s}_{\varepsilon}(z) = \{x \in W^{s}(z):  d(\varphi^{k} (x), 
    \varphi^{k}(z))  \leq \varepsilon \ \forall \ k\geq 0\}.
    \]
A key observation is that the local stable 
set may, under suitable conditions,
have a regular geometrical structure. 
In particular, if \( M \) is a smooth Riemannian manifold 
and \( \varphi \) is a differentiable map, a  typical 
statement  of a  ``Stable Manifold Theorem''
is the following:
\begin{center}\( 
W^{s}_{\varepsilon}(x) \)  \emph{is a smooth submanifold of} \( M \).
\end{center}
This implies in particular that the global stable manifold, which can 
be written as 
\[ 
W^{s}(x) = \bigcup_{n\geq 0} \varphi^{-n}(W^{s}_{\varepsilon}(z)),
\]
is also a  smooth (immersed) submanifold of M; it may however fail to be an
\emph{embedded} manifold (i.e. a manifold in the topology induced from 
\( M \)) due to  the complicated way in which it may twist 
back on itself. 

\subsection{Historical remarks}
   \label{hist}
   As befits such a fundamental result, there exists an enormous 
   literature on the subject, tackling the problem under a number of 
   different conditions. A key 
    idea is that of \emph{hyperbolicity}. In the simplest setting 
   we say that a fixed point \( p \) is hyperbolic if the derivative \( 
   D\varphi_{p} \) has no eigenvalue on the unit circle.  In the 
   analytic, two-dimensional,  
   area preserving case Poincar\'e  proved that the local stable (and 
   unstable) sets are  analytic submanifolds \cite{Poi86s}. Hadamard 
    and Perron  independently developed  more geometric methods 
   allowing them to assume only a \( C^{1} \) 
   smoothness condition \cites{Had01s, Had01, Per28, Per30}; the stable manifold 
   theorem for hyperbolic fixed points is sometimes called the 
   \emph{Hadamard-Perron Theorem}. 
   In \cite{Ste55}, Sternberg used a simple geometric 
   argument, related to Hadamard's technique, 
   to obtain existence and regularity results assuming 
   only \emph{partial} hyperbolicity of the fixed point, i.e. assuming 
   only that the two eigenvalues are real and distinct. 
   Other early work on the subject includes \cites{BogMit61, BogMit63, 
   Hal60, Hal61, Dil60, Dil61} and \cite{Sac67} in which the techniques 
   were generalized to deal with stable manifolds associated to more 
   general compact sets as opposed to just fixed points.  

   In the late 60's and early 70's the theory of stable manifolds became 
   fundamental to the theory of \emph{Uniformly Hyperbolic} dynamics 
   pioneered by Anosov \cite{Ano67} and Smale \cite{Sma67}: there exists
    a \emph{continuous} decomposition 
   \begin{equation*}\label{split}
   T\Lambda = E^{s}\oplus E^{u}
   \end{equation*}
   of the tangent bundle over some set \( \Lambda \) into subbundles on which 
   uniform contraction and exponential estimates hold under the action of 
   the derivative. A straightforward generalizations of this set-up is 
   that  of \emph{partial} or \emph{normal} uniform hyperbolicity
   which allows for the 
   possibility of a neutral subbundle
   \begin{equation*}\label{splitpart}
   T\Lambda = E^{s}\oplus E^{c}\oplus E^{u}.
   \end{equation*}
   This is a significant weakening of the uniform 
   hyperbolicity assumptions as it allows the dynamics tangent to \( E^{c} \) 
   to be quite general. 
  Such situations have been
   systematically and thoroughly investigated using variations and 
   generalizations of the basic methods of Hadamard and Perron
   \cites{Kel67, HirPug69, HirPug70, HirPalPugShu70, Irw70, Fen71, 
   Fen74}, see \cite{HirPugShu77} for a comprehensive treatment. 

   An even more general set-up is based 
   on the \emph{Multiplicative Ergodic Theorem} of Oseledets 
   \cite{Ose68} which says that there always exists 
   a \emph{measurable} decomposition 
   \[ 
   T\Lambda = E^{1}\oplus E^{2} \oplus \cdots \oplus E^{k}
   \]
 with respect to \emph{any} invariant probability measure \( \mu 
   \), such that the \emph{asymptotic exponential growth rate} 
   \[ 
   \lim_{n\to\infty} \frac{1}{n}\log \|D\varphi_{x}^{n}(v)\| = \lambda^{i}
   \]
   is well defined, and for ergodic \( \mu \) even
   independent of \( x \), for all non-zero \( v\in 
   E^{i} \). The condition \( \lambda^{i}\neq 0 \) for all \( i=1,\ldots, 
   k \) is a condition of \emph{non-uniform hyperbolicity} (with respect 
   to the measure \( \mu \)) since it implies that all vectors as 
   asymptotically contracted or expanded at an exponential rate. The 
   non-uniformity comes from the fact that the convergence to the limit 
   is in general highly non-uniform and thus one may have to wait an 
   arbitrarily long time before this exponential behaviour becomes 
   apparent. Pesin \cite{Pes76, Pes77} extended many results of the 
   theory of uniform hyperbolicity concerning stable manifolds 
   to the non-uniform setting. 

   There 
   have also been some recent papers introducing new approaches and 
   focussing on particular subtleties of interest in various contexts 
   \cites{AbbMaj01, Cha01, Cha02}.  We emphasize however that all these 
   results assume \emph{exponential} estimates for the derivative on the 
   relevant subbundles.

\subsection{Very weak hyperbolicity} 

The aim of this paper is to develop some techniques suitable for 
dealing with situations with \emph{very weak} forms of hyperbolicity.
For \( z\in M \) and \( k\geq 1 \)
let
\[
F_{k}(z) = \|(D\varphi^{k}_{z})\| = 
\max_{\|v\|=1}\{\|D\varphi^{k}_{z}(v)\|\}, 
\]
and 
\[
E_{k}(z) =
\|(D\varphi^{k}_{z})^{-1}\|^{-1} = 
\min_{\|v\|=1}\{\|D\varphi^{k}_{z}(v)\|\}.
\]  
These quantities have a  simple geometric interpretation: 
since \( D\varphi_{z}^{k}: T_{z}M \to 
T_{\varphi(z)}M \) is a linear map, it sends circles to ellipses; then
 \( F_{k}(z)\) is precisely half the length of 
major axis of this ellipse and \( E_{k}(z) \) is 
precisely half the length of the minor axis of this ellipse. Then let 
\[
H_{k}(z) = \frac{E_{k}(z)}
{F_{k}(z)}. 
\]
Notice that we always have \( H_{k}(z) \leq 1 \). 
The weakest possible \emph{hyperbolicity condition} one could assume
on the orbit of some point \( x \) is perhaps the condition 
\[
H_{k}(z) < 1
\]
for all \( k \geq 1 \) (or at least all \( k \) sufficiently large), 
equivalent to saying that the image of the unit circle is 
strictly an ellipse or that \( D\varphi^{k}_{z} \) is not 
conformal. 
At the other extreme, perhaps the strongest hyperbolicity condition is
to assume that 
\( H_{k}(z) \to 0 \) exponentially fast in \( k \). This is the case
in the classical  hyperbolic 
setting, both uniform and nonuniform. 
In this paper we prove a
stable manifold theorem essentially under the 
``summable" hyperbolicity condition
\begin{equation*}
\sum H_{k}(z) < \infty.
\end{equation*}
The precise statement of the results 
requires  some additional technical conditions which will be given 
precisely in the next section,  however  
the main idea is that the usual exponential decay of \( H_{k} \) is an 
unnecessarily strong condition. 
Existing
arguments  rely
on a contraction mapping theorem in some suitable space of ``candidate'' 
stable manifolds which yields a fixed 
point (corresponding to the real stable manifold) by the observation that a 
certain sequence is Cauchy and thus converges.  In our 
approach we construct an  \emph{canonical} sequence of 
\emph{finite time local stable manifolds} and use the summability 
condition to show directly that this sequence is Cauchy and 
thus converges to a real stable manifold. 
Also, we make no 
a priori assumptions on the existence of any tangent space 
decomposition.

\subsection{Finite time local stable manifolds}

Our method is based on the key notion of 
\emph{finite time} local stable manifold.  
Let \( k\geq 1 \) and suppose that 
 \( H_{k}(z) < 1 \); then we let 
\( e^{(k)}(z) \) and \( f^{k}(z) \) denote  unit vectors in the 
directions which are 
\emph{most contracted} and \emph{most expanded} respectively by 
\( D\varphi^{k}_{z} \). 
Notice that these directions are 
solutions to
the differential equation \( d\|D\varphi_{z}^{k}(v)\|/d\theta = 0 \) 
    which are given by 
    \begin{equation}\label{contractive directions}
\tan 2\theta  =
\frac{2 (\pfi x1^{k}\pfi y1^{k} +\pfi x2^{k}\pfi y2^{k})}
{(\pfi x1^{k})^2+(\pfi x2^{k})^2 - (\pfi y1^{k})^2 -(\pfi y2^{k})^2}.
\end{equation}
In particular, \( e^{(k)} \) and \( f^{(k)} \) are \emph{orthogonal}
and,  if  \( \varphi^{k} \) is \( C^{2} \), 
\emph{continuously differentiable}  in some 
neighbourhood \( \mathcal N^{(k)}(z) \) in which they are defined. 
Therefore they determine
 two orthogonal foliations \( \mathcal E^{(k)}\) and \( \mathcal
 F^{(k)} \) defined by the integral curves of the unit vector fields
\(e^{(k)}(x) \) and \( f^{(k)}(x) \) respectively. We let 
 \( \mathcal E^{(k)}(z) \) and \( \mathcal F^{(k)}(z) \) 
 denote  the corresponding leaves through the point \( z \). 
These are the natural finite time versions of the local stable and 
unstable manifolds of the point \( z \) since they are, in some 
sense, the most contracted and most expanded curves through \( z \) 
in  \( \mathcal N^{(k)}(z) \). Notice that they are uniquely defined 
locally.  We will show 
that under suitable conditions the finite time local stable 
manifolds converge to  
real local stable manifold.

The idea of constructing finite time 
local stable manifolds is not new. 
In the context of Dynamical Systems, 
as far as we know it was first introduced in 
 \cite{BenCar91} and developed further in several papers 
 including \cites{MorVia93, BenYou93, LuzVia2, HolLuz, WanYou01} 
 in which systems satisfying some nonuniform hyperbolicity are 
 considered. 
All these papers deal with families of systems in which, initially, 
hyperbolicity cannot be guaranteed for all time for all parameters. A 
delicate parameter-exclusion argument 
requires  information about the geometrical structure of stable and 
unstable leaves based only on a finite number of iterations and thus 
the notion of finite time manifolds as given above is very natural. 
We emphasise however that in these papers
the construction is heavily embedded in the global 
argument and no particular emphasis is placed on this method as an 
algorithm for the construction of real local stable manifolds \emph{per 
se}. Moreover the decay rate of \( H_{k} \) there is exponential 
and the specific properties of the systems 
(such as the small determinant and various other hyperbolicity 
and distortion conditions) are heavily 
used, obscuring the precise conditions required for the argument to work.

One  aim of this paper is to 
clarify the setting and assumptions 
required for the construction to work and to show that the 
main ideas can essentially be turned into a 
fully fledged alternative approach to 
theory of stable manifolds. Moreover we show that the argument 
goes through under much weaker conditions than those 
which hold in the papers cited above.

\subsection{Main Results}
\label{results}

We shall consider dynamical systems given by maps 
\[
\varphi: M \to M 
\]
where
\( M \) is a two-dimensional
Riemannian manifold with Riemannian metric \( d \). 
The situation we have in mind  is that of 
a piecewise \( C^{2} \) diffeomorphism
with singularities: 
there exists a set \( \mathcal 
S \) of zero measure such that \( \varphi \) is a \( C^{2} \) local 
diffeomorphism on \( M\setminus\mathcal S \). 
The map \( \varphi \) may be discontinuous on \( \mathcal S \) and/or 
the first and second derivatives may become unbounded near \( \mathcal 
S \).  The precise  assumptions will be local and will be formulated 
below. 
First of all we introduce some notation. For  \( x \in M \) let 
\[
P_k(x)=\|D\varphi_{\varphi^kx}\|,\quad
Q_k(x)=\|(D\varphi_{\varphi^kx})^{-1}\|,\quad
\widetilde{P}_k(x)=\|D^2\varphi_{\varphi^kx}\|,
\]
and
\[
\dd_k(x)=|\det D\varphi_{\varphi^kx}|,\quad
\tilde\dd_k(x)=\|D(\det D\varphi_{\varphi^kx})\| 
\]
Notice that all of these quantities 
 depend only on the derivatives 
of \( \varphi \) at the point \( \varphi^{k}(x) \). If \( \varphi \) 
is globally a \( C^{2} \) diffeomorphism then they are all uniformly 
bounded above and below and play no essential role in the result. 
On the other hand, if the contraction and/or expansion is unbounded 
near the singularity set \( \mathcal S \) some control of the 
recurrence is implicitly given by some conditions which we impose on 
these quantities. 
We shall also use the notation
\[   F_{j, k}(x) = \|D\varphi^{k-j-1}_{\varphi^{j+1}(x)}\|.  \]
We now give a generalization of the notion of local stable manifold. 
For a sequence
\[ 
\underline\varepsilon=\{\varepsilon_j\}_{j=0}^{\infty} 
\]
with 
\(
\varepsilon_{j}\geq \varepsilon_{j+1}> 0
\)
for all  \( j \geq 0 \), 
we  let 
\[
\mathcal N^{(k)} = \mathcal N^{(k)}_{\underline\varepsilon}(z) 
=\{\tilde{z}\in M:\,\|\varphi^{j}(\tilde{z})-\varphi^{j}(z)\|\leq
\varepsilon_j,\forall j\leq k-1\}. 
\]
This defines a nested sequence of neighbourhoods of the 
point \( z \). 
We shall always suppose that for all \( k\geq j \geq 1 \)  the 
restriction \( \varphi^{j}|
\mathcal N^{(k)} \) is a \( C^{2} \) diffeomorphism   onto its
image. In the presence of singularities this may impose a strong 
condition on the sequence \( \underline \varepsilon \) whose terms may 
be required to decrease very quickly. 
We then let 
\( \{
p_{k}, 
q_{k}, \tilde{p}_{k} \}_{k=1}^{\infty}
\)
be uniform upper bounds for the values of \( P_{k}, Q_{k}, \tilde P_{k} 
\) respectively in \( \mathcal N^{(k)}(z) \): 
\[ 
p_{k}=\max_{x\in\mathcal N^{(k)}}P_{k}(x), \quad 
q_{k}=\max_{x\in\mathcal N^{(k)}}Q_{k}(x), \quad 
\tilde p_{k}=\max_{x\in\mathcal N^{(k)}}\tilde P_{k}(x). 
\]
These values may be unbounded.  
Then let 
 \(
 \{ \gamma_{k},  \gamma^{*}_{k}, 
 \delta_{k} \}_{k=1}^{\infty}
\) be given by 
\[ 
 \gamma_{k}= \max_{x\in\mathcal N^{(k)}}
 \{H_{k}\},\quad \gamma^{*}_{k}= \max_{x\in\mathcal N^{(k)}}\{E_{k}\},
\]
and 
\[ 
\delta_{k}= 
\max_{x\in\mathcal N^{(k)}} 
\left\{
\frac{E_k}{F^{2}_{k}}\sum_{j=0}^{k-1} \tilde{P}_{j}F_{j, k}
F^{2}_{j} + 
\frac{E_{k}}{F_{k}} \sum_{j=0}^{k-1}
\dd^{-1}_{j}\tilde{\dd}_{j}F_{j} 
\right\}.
\]
We are now ready to state our two hyperbolicity conditions. The first 
is a hyperbolicity condition
    \begin{equation}\tag{\( * \)}
\sum_{k=1}^{\infty} 
p_k q_k\gamma_{k+1}+ 
\tilde{p}_k q^{5}_{k} p_k^3\gamma^{*}_{k+1} +
p_k^{5}q_k^5\delta_{k} + p_k^{2}q_k^2\delta_{k+1}< \infty.
\end{equation}
Notice that if the norm of the derivative is bounded, such as 
in the absence of singularities, this
reduces to the more ``user-friendly'' condition 
\[ 
\sum \gamma_{k}+\gamma^{*}_{k}+\delta_{k} < \infty. 
\]
The summability of \( \{\gamma^{*}_{k}\} \) is not particularly
crucial and is really only used to ensure that some minimal
contraction is present, so that the presence of a contracting stable
manifolds makes sense. The summability of \( \{\gamma_{k}\} \) is
simply the ``summable hyperbolicity'' assumptions stated above. The
summability of \( \{\delta_{k}\} \) is a quite important technical
assumption related to the ``monotonicity'' of the estimates in \( k
\), it is not overly intuitive but it is easily verified 
in standard situations such as in the uniformly hyperbolic setting. 
Taking advantage of condition \( (*) \) we define
\[
 k_{0} =  \min\{j: p_k q_k\gamma_{k+1}< 1/2,\,
 \forall\, k\geq j-1\}\} < \infty
\]
and the sequence 
\[
\tilde\gamma_{k} = \gamma^{*}_k+
2 \max_{x\in\mathcal N^{(k)}} 
\{F_{k}\}\sum_{i=k}^{\infty}p_iq_i\gamma_{i+1} < \infty.
\]
Our second assumption is that there exists some constant 
\( \Gamma >0 \) such that 
\begin{equation*}\tag{\( ** \)}
    \tilde\gamma_{j} + 4 \max_{\mathcal N^{(j)}} \{\| F_{j}\|\}
    p_{k}q_{k}\gamma_{k+1} < 
\Gamma \varepsilon_{j}
\end{equation*}
for all   \( k\geq k_{0}  \) and  \( j\leq k \). 

This is not a particularly intuitive condition but 
thinking of it in the simplest setting can be useful. 
Supposing for example that we are in a uniformly hyperbolic situation 
and that all derivatives are bounded, we have that
the left hand side is \(  \lesssim E_{k} \) which specifies that in
some sense, the images of the neighbourhoods of \( z \) under consideration 
should not shrink to fast relative to the contraction in these
neighbourhoods.
We now state our main result. 

\begin{main theorem}
    Let \( z\in M \) and suppose that there exists a sequence 
    \( \underline \varepsilon \) such that \( \varphi^{k} \) 
    restricted to \( \mathcal N^{(k)} \) is a \( C^{2} \)
   diffeomorphism onto its image for all \( k\geq 1 \), and suppose 
   also that  conditions \( (*) \) and \( (**) \) hold. 
Then there exists \( 
    \varepsilon>0 \) and a
\( C^{1+\text{Lip}} \) embedded one dimensional submanifold 
\( \mathcal{E}^{\infty}(z)  \) of \( M 
\) containing \( z \) such that \( |\mathcal E^{(\infty)}(z)| \geq \varepsilon 
\) and such that 
there exists a constant \( C>0 \) such that 
\( \forall z, z' \in \mathcal E^{\infty}(z) \ \forall \ k\geq k_0 \) we have
\[ 
\textstyle{|\varphi^{k}(z)- \varphi^{k}(z')| \leq C \tilde\gamma_{k} |z -z'|.} 
\]
In particular if $\tilde\gamma_k\to 0$  then \(  |\varphi^{k}(z)- \varphi^{k}(z')| \to 0\) as 
\( k\to\infty \) and therefore 
\[
\mathcal E^{(\infty)}(z) \subseteq W^{s}_{\varepsilon} (z). 
\]
Moreover if
\( F_{k}\to\infty \) uniformly in \( k \), 
then 
\[ \mathcal E^{(\infty)}(z) = 
\bigcap_{k\geq k_0}\mathcal{N}^{(k)}(z).\]
   \end{main theorem}

   We divide the proof into several sections. In \ref{section_notation} 
   we introduce some useful notation. In \ref{s:dist} we prove a 
   technical estimate which shows that the summability condition on \( 
   \delta_{k} \) implies some uniform distortion bounds on the \( 
   \mathcal N^{(k)} \). In \ref{point_theory} we study the convergence of 
   pointwise contracting directions and in \ref{sec-regularity1} 
   we use these to study the convergence of the local finite time stable 
   manifolds. In \ref{strategygeom} we show that the limit curve has 
   positive length. This is not directly implied by the preceding 
   convergence estimates which give convergence of the leaves on 
   whichever domain they are defined. Here we need to make sure that 
   such a domain of definition (i.e. length) of the leaves can be 
   chosen uniformly. Thus we have to worry about the shrinking of the 
   sets \( \mathcal N^{(k)}(z)  \). Condition \( (**) \) is used 
   crucially in this section. We remark that the  lower bound \( \varepsilon \)
   for the length of the local stable manifold is determined in this
   section. 
   In \ref{sec-regularity2} we show that the 
   limit curve is smooth and in \ref{sec-contraction} that it 
   ``contracts'' and is therefore indeed part of the local stable 
   manifold. Finally, in \ref{uniqueness} we discuss uniqueness issues.

\section{Hyperbolic fixed points}
 \label{fixed_point_case}

 As an application of our abstract theorem, 
 we consider the simplest case of a hyperbolic fixed
 point. The result is of course already well-known in this context, but
 we show  that our conditions are easy to check and that it
 therefore follows almost immediately from our general result. 

 Let \( M \) be a two-dimensional
 Riemannian manifold with Riemannian metric \( d \), and let \(
 \varphi: M \to M \) be a   \( C^{2} \) diffeomorphism. 
 Suppose that \( p\in M \) is a fixed point. The local stable manifold of $p$ is the set
 \( W^{s}(p) \)  of points which 
 remain in a fixed neighbourhood of \( p \) for all forward 
 iterations: for \( \eta > 0 \) and \( k\geq 1 \) let 
    \[ 
    \mathcal N^{(k)}_{\eta}=\{x: d(\varphi^{j}(x), p) \leq 
    \eta \ \forall \ 0\leq j \leq k-1\}
 \]
 and
 \[
    \mathcal N^{(\infty)}_{\eta}= \bigcap_{k\geq 1} \mathcal 
     N^{(k)}
    \]
	 For  \( \eta > 0 \), 
 we define the \emph{local stable set} of  \( 
 p \) by
 \[
     W^{s}_{\eta}(p) = W^{s}(p) \cap\mathcal 
     N^{(\infty)}_{\eta}(p)
 \]

 In this section we shall focus on the simplest setting 
 of a hyperbolic fixed point. 
     We recall that the fixed point \( p \) is hyperbolic if the 
     derivative \( D\varphi_{p} \) has no eigenvalues on the unit circle. 
     
  \begin{theorem*}
      Let \( \varphi: M \to M \) be a \( C^{2} \) diffeomorphism of a 
      Riemannian surface and suppose that \( p \) is a hyperbolic 
      fixed point with eigenvalues 
 \( 0<| \lambda_{s}|< 1 < | \lambda_{u}| \). Given $\eta>0$,
      there exists a constant 
      \( \varepsilon(\eta) > 0 \) such that the following properties hold: 
      \begin{enumerate}
 \item \( W^{s}_{\eta} (p) \) is \( C^{1+Lip} \) 
	  one-dimensional submanifold of M 
	 tangent to \( E^{s}_{p} \);
 \item \( |W^{s}_{\eta} (p)|\geq \varepsilon \) 
	 on either side of \( p \);
 \item \( W^{s}_{\eta} (p) \) contracts at an 
	 exponential rate.
	 \item 
     \[ 
 W^{s}_{\eta}(p) = \bigcap_{k\geq 0} \mathcal N^{(k)}_{\eta}(p). 
 \]	
 \end{enumerate}
      \end{theorem*}

 \begin{proof}
 To prove this result, it suffices to verify the hyperbolicity conditions 
 stated in section \ref{results}.
 First of all , since $\varphi$ is a $C^2$ diffeomorphism, all the first and second partial derivatives
 are continuous and bounded. Hence for all $k\geq 0$, there is a uniform constant $K>0$ such that 
 $$p_k,q_k,\tilde{p}_k,\dd_k\tilde{\dd}_k\leq K.$$ To estimate expansion and
 contraction rates in $\mathcal N^{(k)}_{\eta}$ we have the following lemma:

 \begin{lemma}\label{regular growth}
    There exists a constant \( K>0 \) such that 
 for all \( \delta  > 0 \) there exists 
 \( \eta (\delta) > 0 \) 
 such that for 
 all \( k\geq j\geq 0 \) and all \( x\in \mathcal N^{(k)}_{\eta} 
 \) we have
 \begin{equation}\label{eguh}
 K ( \lambda_{u}+ \delta)^{j} \geq F_{j} \geq 
 ( \lambda_{u}- \delta)^{j} \geq 
 (\lambda_{s}+ \delta)^{j} 
 \geq  E_{j} \geq 
 K^{-1} (\lambda_{s}- \delta)^{j}
 \end{equation}
 and  
 \begin{equation}\label{eguh1}
 \sum_{j=0}^{k-1} F_{j} \leq K F_{k}; \quad
 F_{j} F_{j, k} \leq K F_{k}; \quad\text{ and } \quad
 \sum_{i =j}^{\infty} H_{i} \leq K H_{j} 
 \end{equation}
 In particular 
 \begin{equation}\label{DUH}
 \|D^{2}\varphi^{k}_{x}\| \leq K F^{2}_{k}; 
 \quad \text{ and } \quad
 \|D(\det
 D\varphi^{k}_{x})\|
 \leq K E_{k} F_{k}^{2}.  
 \end{equation}

 \end{lemma}

 \begin{proof}
 The estimates in \eqref{eguh} and (\ref{eguh1}) 
     follow from standard estimates in the theory of uniform 
     hyperbolicity. We refer to \cite{KatHas94} for details and proofs. 
     The estimates in \eqref{DUH} then follow from substituting 
     \eqref{eguh1} into the estimates of Lemma \ref{second derivative}.
 \end{proof}

 Next we verify hyperbolicity conditions $(*)$ and $(**)$. We estimate $\gamma_k,\tilde\gamma_k,\gamma^{*}_k$ 
 and $\delta_k$ for each $k\geq 0$. For $\gamma_k$ and $\gamma^{*}_k$ we have
 \begin{gather*}
 \gamma_k=\max_{x\in\mathcal N^{(k)}}
  \{H_{k}\}\leq K\frac{(\lambda_{s}+ \delta)^{k}}{(\lambda_{u}-\delta)^{k}},\\
 \gamma^{*}_k=\max_{x\in\mathcal N^{(k)}}
  \{E_{k}\}\leq K(\lambda_{s}+ \delta)^{k},
 \end{gather*}
 while for $\tilde\gamma_k$ we obtain
 \begin{equation*}
 \begin{split}
 \tilde\gamma_k &= \gamma_k+
 \max_{x\in\mathcal N^{(k)}} 
 \{2F_{k}\}\sum_{i=k}^{\tilde k-1}p_iq_i\gamma_{i+1}\\ 
 &\leq 2(\lambda_{s}+ \delta)^{k}+K\biggl[\frac{(\lambda_{u}+\delta)(\lambda_{s}+\delta)}
 {(\lambda_{u}-\delta)}\biggr]^k\leq 
 K(\lambda_{s}+\tilde\delta)^k,
 \end{split}
 \end{equation*}
 where $\tilde\delta$ can be made small with $\delta$ small.
 To estimate $\delta_k$, we just use Lemma \ref{regular growth} above to conclude that
 \begin{equation*}
 \begin{split}
 \delta_{k} &= 
 \max_{x\in\mathcal N^{(\tilde k)}} 
 \left\{
 \frac{E_k}{F^{2}_{k}}\sum_{j=0}^{k-1} \tilde{P}_{j}F_{j, k}
 F^{2}_{j} + 
 \frac{E_{k}}{F_{k}} \sum_{j=0}^{k-1}
 \dd^{-1}_{j}\tilde{\dd}_{j}F_{j} 
 \right\}\\
 &\leq  K(\lambda_{s}+\tilde\delta)^k.
 \end{split}
 \end{equation*}
 In the estimates above, the constant $K$ is uniform and depends only
 on $\lambda_s, \lambda_u$ and the bounds for the partial derivatives of $\varphi$.

 Condition $(*)$ is now immediate, since for $\delta,\tilde\delta$ sufficiently small, the constants 
 $\gamma_k,\gamma^{*}_k,\tilde\gamma_k$ and $\delta_k$ all decay exponentially fast. In particular there 
 exists a constant $L>0$ such that $\mathrm{Lip}(e^{(k)})\leq L$ inside each $\mathcal{N}^{(k)}$.

 Let $k_0$ be the constant defined in section \ref{results}. To verify condition $(**)$, 
 we just need to show that there is a $\Gamma>0$ such that $\forall k\geq k_0$ we have:
 \begin{equation}\label{geom_fixedpt}
 K\frac{(\lambda_{u}+\delta)^{j}(\lambda_{s}+ \delta)^{k+1}}
 {(\lambda_{u}-\delta)^{k+1}}+K(\lambda_{s}+\tilde\delta)^k<\Gamma\eta,\qquad\forall j\leq k+1.
 \end{equation}
The existence of $\Gamma$ follows immediately if we
choose $\delta$ sufficiently small so that 
$$(\lambda_{u}+\delta)(\lambda_{s}+ \delta)(\lambda_{u}-\delta)^{-1}<1,\quad\textrm{and}
\quad\lambda_{s}+\tilde\delta<1.$$
The conclusions of the theorem now follow.
In particular, the length $\varepsilon$ of the limiting leaf $\mathcal{E}^{(\infty)}$
is determined by equation (\ref{epsilon1}) in section \ref{geom_section}.
\end{proof}

\section{Finite time local stable manifolds}

In this section we prove some estimates concerning the relationships
between finite time local stable manifolds of different orders. In
particular we prove that they form a Cauchy sequence of smooth curves.
Throughout this and the following 
 section we work under the assumptions of our main 
 theorem. In particular we consider the orbit of a point \( z \) and 
 are given a sequence of neighbourhoods \( \mathcal N^{(k)}=\mathcal 
 N^{(k)}(z) \) in which most contractive and most expanding directions 
 are defined and thus, in particular, in which the finite time local stable 
 manifolds \( \mathcal E^{(k)}(z) \) are defined. The key problem 
 therefore is to show that these finite time local stable manifolds 
 converge, that they converge to a smooth curve, and that this curve 
 has non-zero length !

\subsection{Notation}
\label{section_notation}
We shall use \( K \) to denote a generic constant which is
allowed to depend only on the diffeomorphism \( \varphi \).
For any    \( j\geq 1 \) we  let 
    \[
    e^{(k)}_{j}(x) =
    D\varphi^{j}_{x}(e^{(k)}(x)) 
    \quad \text{ and } \quad
    f^{(k)}_{j}(x) =
    D\varphi^{j}_{x}(f^{(k)}(x)) 
    \]
denote the images of the most contracting and most expanding vectors. 
To simplify the formulation of angle estimates we introduce the 
variable \( \theta \) to define the position of the vectors. We write
\begin{gather*}
e^{(n)}=(\cos\t^{(n)},\,\sin\t^{(n)}),\quad
f^{(n)}=(-\sin\t^{(n)},\,\cos\t^{(n)}).\\
e^{(n)}_{n}=E_n (\cos\t^{(n)}_{n},\,\sin\t^{(n)}_{n}),\quad
f^{(n)}_{n}=F_n (-\sin\t^{(n)}_{n},\,\cos\t^{(n)}_{n}).
\end{gather*}
Finally, we let 
\[\phi^{(k)}=\measuredangle(e^{(k)},e^{(k+1)})
\text{ and } 
\phi^{(k)}_j=\measuredangle(e^{(k)}_j,e^{(k+1)}_j). 
\]
We also identify any vector $v$ with $-v$, or equivalently we identify 
an angle $\theta$ with the angle $\theta+\pi$. 
Important parts of the proof depend on estimating the derivative of 
various of these quantities with respect to the base point \( x \). We 
shall write \( D\phi^{(k)}, De^{(k)}, D\theta^{(n)}_{j}, \ldots  \) to denote 
the derivatives with respect to the base point \( x \). 
To simplify the notation we let 
\begin{equation}
\Xi_{k}(x) := 
\frac{P_k(x)Q_k(x) H_{k+1}(x)}{(1-P_k(x)Q_k(x) H_{k+1}(x))} 
\leq \frac{p_kq_k\gamma_{k+1}}{(1-p_kq_k\gamma_{k+1})} :=\xi_k
\end{equation}
Also, all statements 
hold uniformly for all \( x\in \mathcal N^{(k)} \). 
 
\subsection{Distortion}
\label{s:dist}

The following distortion estimates follow from completely general 
calculations which do not depend on any hyperbolicity assumptions. 
The definition of \( \delta_{k} \) is motivated by these estimates 
which will be used extensively in section \ref{point_theory}. 

\begin{lemma}\label{second derivative}
For all \( k\geq 1 \) and and all \( x \) such that \( \varphi^{k} \) 
is \( C^{2} \) at \( x \), we 
have 
\begin{equation*}\tag{D1}\label{D1}
H_k\frac{ \|D^{2}\varphi^{k} \|}{\|D\varphi^{k}\|} \leq 
\frac{E_k}{F^{2}_{k}}\sum_{j=0}^{k-1} \tilde{P}_{j}F_{j, k}
F^{2}_{j} \ \ (\leq \delta_{k})
\end{equation*}
and
\begin{equation*}\tag{D2}\label{D2}
\frac{\| D(\det D\varphi^{k}_{z})\|}
{\|D\varphi^{k}\|^{2}} \leq  \frac{E_{k}}{F_{k}} \sum_{j=0}^{k-1}
\dd^{-1}_{j}\tilde{\dd}_{j}F_{j} \ \ (\leq \delta_{k})
\end{equation*}
   \end{lemma}

\begin{proof}
       Let \( A_{j}= D\varphi_{\varphi^{j}z} \) and let \( A^{(k)}=
       A_{k-1}A_{k-2}\dots A_{1} A_{0} \).  Let \( D A_{j}\) denote
       differentiation of \( A_{j} \) with respect to the space variables. 
       By the product rule for differentiation we have
\begin{equation}\label{productrule}
   \begin{split}
D^{2}\varphi^{k}_{z} &= DA^{(k)} 
= D(A_{k-1}A_{k-2}\dots
   A_{1}A_{0}
   ) \\ &= \sum_{j=0}^{k-1} A_{k-1}\dots A_{j+1}(DA_{j}) A_{j-1}\dots A_{0}.
   \end{split}
\end{equation}
Taking norms on both sides of 
\eqref{productrule} and using the fact that 
\( A_{k-1}\dots A_{j+1} = D\varphi^{k-j-1}_{\varphi^{j+1}z} \), \(
A_{j-1}\dots A_{0} = D\varphi^{j}_{z} \) and, by the chain rule, \(
DA_{j} = D (D\varphi_{\varphi^{j}z}) = D^{2}\varphi_{\varphi^{j}z}
D\varphi^{j}_{z} \), we get
\[
\|D^{2}\varphi^{k}_{z}\| \leq \sum_{j=0}^{k-1}
\|D\varphi^{k-j-1}_{\varphi^{j+1}z}\| \cdot \|
D^{2}\varphi_{\varphi^{j}z}\| \cdot \|D\varphi^{j}_{z} \|^{2} \leq 
\sum_{j=0}^{k-1} \| D^{2}\varphi_{\varphi^{j}z}\| F_{j,k} F_{j}^{2}.
\]
The inequality (D1) now follows. For (D2) 
we argue along similar lines, this time letting \(
A_{j}= \det D\varphi_{\varphi^{j}z} \).  Then we have, as in
\eqref{productrule} above,
$$
D(\det\varphi^{k}_{z}) = DA^{(k)} = \sum_{j=0}^{k-1} A_{k-1} \dots
A_{j+1}(DA_{j}) A_{j-1}\dots A_{0}.$$
Moreover we have that
\( A_{k-1}\dots A_{j+1} = \det D\varphi^{k-j-1}_{\varphi^{j+1}z} \),
\( A_{j-1}\dots A_{0} = \det D\varphi^{j}_{z}, \) and by the chain
rule, also: $$ DA_{j} = D(\det D\varphi_{\varphi^{j}z}) = (D \det
D\varphi_{\varphi^{j}z}) D\varphi^{j}_{z}.$$ This gives
\begin{equation}\label{det}
   D(\det\varphi^{k}_{z}) = \sum_{j=0}^{k-1} (\det
   D\varphi^{k-j-1}_{\varphi^{j+1}z}) (D \det D\varphi_{\varphi^{j}z})
   (\det D\varphi^{j}_{z} )(D\varphi^{j}_{z}).
\end{equation}
By the multiplicative property of the determinant we have
the equality: $$
(\det D\varphi^{k-j-1}_{\varphi^{j+1}z}) (\det D\varphi^{j}_{z} ) =
\det D\varphi^{k}_{z}/\det D\varphi_{\varphi^{j}{z}}.$$ Thus, taking
norms on both sides of \eqref{det} gives
\[ 
\| D(\det D\varphi^{k}_{z})\| \leq \ |\det D\varphi^{k}_{z} |
\sum_{j=0}^{k-1} \frac{\|D(\det D\varphi_{\varphi^{j}(z)})\|}{|\det
D\varphi^{\varphi^{j}(z)}|} F_{j}
\]
The inequality (D2)  now follows
from the fact that \( \det D\varphi^{k} = E_{k} F_{k}\).
\end{proof}

\subsection{Pointwise convergence}
\label{point_theory}

In this section we prove two key lemmas showing that both the angle 
\( \phi^{(k)} \) (Lemma \ref{angleconvergence}) 
between consecutive most contracted directions and 
the norm of its spatial derivative \( D\phi^{k} \) (Lemma 
\ref{derivativeconvergence}) can be bounded in terms of the 
hyperbolicity. In particular, from the summability condition \( (*) 
\), we obtain also that the norm \( 
\|De^{(k)}\| \) of the spatial 
derivative  of the contractive directions is uniformly bounded in \( k 
\). 

\begin{lemma} \label{angleconvergence}
    For all \( k\geq k_{0} \) and \( x\in\mathcal N^{(k)} \) we 
    have 
 \begin{equation}\label{angle}
 |\phi^{(k)}| \leq |\tan\phi^{(k)}| \leq\frac{P_kQ_k H_{k+1}}{
1-P_kQ_k H_{k+1}} \ \ (\leq \xi_{k}).
\end{equation}
Moreover, for all \( k\geq j\geq  k_{0} \) we have 
\begin{equation}
\label{contractionclaim} 
\|e^{(k)}_{j}(x)\|\leq E_{j}(z)+F_{j}(z)\sum_{i=j}^{k-1}\phi^{(i)}(z)
\ \ 
(\leq \tilde{\gamma}_j).
\end{equation}
\end{lemma}
Notice that the estimate in \eqref{contractionclaim} 
gives an upper bound for the contraction which depends only on \( j \)
and not on \( k \). 
\begin{proof}
    We claim first of all that for all $k\geq k_0$ we have 
\begin{equation}\label{half}
    \|e^{(k)}_{k+1}\|/F_{k+1}\leq P_kQ_k H_{k+1}\leq 1/2. 
\end{equation}
To see this observe that
\(
E_k \leq \|e^{(k+1)}_{k}\|\leq
\|D\varphi^{-1}_{z_k}e^{(k+1)}_{k+1}\|\leq Q_k E_{k+1}\) ,
\( E_{k+1} \leq \|e^{(k)}_{k+1}\|=\|D\varphi_{z_k}e^{(k)}_{k}\|
\leq P_k E_k \),
\( F_k= \|D\varphi^{-1}_{z_k}f^{(k)}_{k+1}\|\leq Q_k F_{k+1} \),
\( F_{k+1} =\|D\varphi_{z_k}f^{(k+1)}_{k}\|\leq P_k F_k \). 
Moreover \( H_{k+1}/H_{k}=(E_{k+1}/F_{k+1})/(E_{k}/F_k) \). 
Combining these inequalities gives 
\begin{equation}\label{minidistortion}
E_{k+1}/E_k\in[Q^{-1}_k,P_k], \quad F_{k+1}/F_k\in[Q^{-1}_k,P_k], 
\quad H_{k+1}/H_k\in[(P_kQ_k)^{-1},P_kQ_k].
\end{equation}
Therefore, writing
write $e^{(k)}_{k+1}=D\varphi(z_k)e^{(k)}_{k}$ and applying the
Cauchy-Schwarz inequality gives 
   \( \|e^{(k)}_{k+1}\|\leq E_{k+1}P_kQ_k \) which immediately implies 
   the first inequality of \eqref{half}. The second inequality follow 
   simply by our choice of \( k_{0} \). 

Now write $e^{(k)}=\eta e^{(k+1)}+\varphi f^{(k+1)}$ where
$\eta^2+\varphi^2=1$ by normalization.
Linearity implies that $e^{(k)}_{k+1}=
\eta e^{(k+1)}_{k+1}+\varphi f^{(k+1)}_{k+1}$ and orthogonality implies that
\(
\|e^{(k)}_{k+1}\|^2 =\eta^2\|e^{(k+1)}_{k+1}\|^2+\varphi^2
\|f^{(k+1)}_{k+1}\|^2
 =\eta^2 E^{2}_{k+1}+\varphi^2 F^{2}_{k+1}
\)
where $E_k=\|e^{(k)}_k\|,\,F_k=\|f^{(k)}_k\|$.
Since $\phi^{(k)}=\tan^{-1}(\varphi/\eta)$ we get
$$
|\tan\phi^{(k)}|=
\left(\frac{\|e^{(k)}_{k+1}\|^2-E^{2}_{k+1}}
{F^{2}_{k+1}-\|e^{(k)}_{k+1}\|^2}\right)^{\frac{1}{2}} \leq
\frac{\|e^{(k)}_{k+1}\|/F_{k+1}}
{\left(1-\|e^{(k)}_{k+1}\|^2/F^{2}_{k+1}\right)^{\frac{1}{2}}}
\leq
\frac{\|e^{(k)}_{k+1}\|/F_{k+1}}
{1-\|e^{(k)}_{k+1}\|/F_{k+1}}.
$$
In the last inequality we have used 
\( \|e^{(k)}_{k+1}\| < F_{k+1} \) from \eqref{half} 
and then applied the inequality
$\sqrt{(1-x^2)}>1-x$, valid for $x\in(0,1)$. 
Using \eqref{half} again
completes the proof of the first statement in the lemma. 

To prove \eqref{contractionclaim} 
we write 
$e^{(k)}_j=e^{(j)}_{j}+\sum_{i=j}^{k-1}(e^{(i+1)}_{j}-e^{(i)}_{j})$. 
The first term is equal to \( E_{j}(x) \) by definition. For the
second we have, by linearity, 
\(
\|e^{(i+1)}_{j}-e^{(i)}_{j}\| = \| F_{j}(x) (e^{(i+1)} - e^{(i)}) \|
\leq \|F_{j}(x)\| \ |\phi ^{(i)}|.
\) By \eqref{angle} and the definition of \( \tilde\gamma_{j} \) we
get \( \| e^{(k)}_j\|\leq \tilde\gamma_{j} \). 
\end{proof}

\begin{lemma}\label{derivativeconvergence}
For all \( k\geq k_{0}  \) and \( x\in \mathcal N^{(k)} \) we have
\[
\|D\phi^{(k)}\|\leq 1597(p_kq_k)^2\delta_{k+1}+40(p_kq_k)^5\delta_{k}
+40(p_kq_k)^3q_{k}^2\tilde{p}_k\gamma^{*}_{k+1}.
\]
In particular, there exists a constant L independent of $k$ such that 
\[
\|De^{(k)}\| \leq
\sum_{j\leq k} \|D\phi^{(j)}\| \leq L.
\]
\end{lemma}
\begin{proof}
    Since \(D\varphi^{k} \) is a linear map, we have
     \begin{equation} \label{hyp0}
     \tan \phi^{(k)} = H_{k+1} \tan \phi^{(k)}_{k+1}.
     \end{equation}
Differentiating \eqref{hyp0} on
both sides and taking norms we have
     \begin{equation}\label{hyp3}
     \begin{split}
     \|D\phi^{(k)}\| &\leq  \| H_{k+1} \cdot D (\tan
     \phi_{k+1}^{(k)})\| + \| D H_{k} \cdot \tan \phi^{(k)}_{k+1}\| \\
     &\leq 
     \|H_{k+1} (1+\tan^{2}\phi^{(k)}_{k+1}) D
     \phi_{k+1}^{(k)}\| + \|D H_{k+1} \cdot \tan \phi^{(k)}_{k+1} \|
     \end{split}
     \end{equation}

In the next two sublemmas we obtain upper bounds for 
\( \|D\phi_{k+1}^{(k)}\| \) and \( \|D H_{k+1} \| \) respectively and 
then substitute these bounds into \eqref{hyp3}. 

\begin{sublemma}\label{derbd}
\[
\|D\phi^{(k)}_{k+1}\|\leq \frac{2048\delta_{k+1}}{9H_{k+1}}+
8(P_kQ_k)^2\frac{\delta_{k}}{H_{k}}  
+8Q^{2}_kP_k\tilde{P}_k F_{k}. 
\] 
  \end{sublemma}
\begin{proof}
Writing \( \phi^{(k)}_{k+1}=\theta^{(k+1)}_{k+1}-\theta^{(k)}_{k+1}
    \) we have  
\[ \|D\phi^{(k)}_{k+1}\|=\|D\theta^{(k+1)}_{k+1}
    -D\theta^{(k)}_{k+1}\| \leq \|D\theta^{(k+1)}_{k+1}\|
    + \|D\theta^{(k)}_{k+1}\|
    \]
    Our strategy therefore is to obtain estimates for the terms on
    the right hand side.
First of all we write
\begin{equation*}
D\varphi^{n}(z)=
\begin{pmatrix}
A_n & B_n \\
C_n & D_n
\end{pmatrix}
\end{equation*}
where $A_n,B_n,C_n$ and $D_n$ are the matrix entries for the derivative
$D\varphi^{n}$ evaluated at $z$.
Since $\{e^{(n)}(z),f^{(n)}(z)\}$ correspond
to (resp.) maximal contracting and expanding vectors under 
$D\varphi^n(z)$, i.e. solutions of the differential equation
\begin{equation*}
\frac{d}{d\t}\left\| D\varphi^{n}_{z}\;
\begin{pmatrix}
\cos\t \\
\sin\t
\end{pmatrix}
\right\|=0
\end{equation*}
By
solving the differential equation above in $\theta$ we get
   \[ 
   \tan 2\t^{(k)} =
   \frac{2(A_k B_k+C_k D_k)}{A^{2}_{k}+C^{2}_{k}-B^{2}_{k}-D^{2}_{k}}
   =\frac{2\A_k}{\B_k}
\]
and  by
solving a similar one for the inverse map $D\varphi^{-n}$ we get 
\[
   \tan 2\t^{(k)}_{k} =
\frac{2(B_k D_k+A_k C_k)}{D^{2}_{k}+C^{2}_{k}-A^{2}_{k}-B^{2}_{k}}
=-\frac{2{\mathcal{C}}_k}{{\mathcal{D}}_k}
   \]   
Notice the use of \( \mathcal A_{k}, \mathcal B_{k}, \mathcal C_{k},
\mathcal D_{k} 
\)  as a shorthand notation for the expression in the quotients.
Now $e^{(k)}_{k},\,f^{(k)}_{k}$ are respectively maximally expanding and
contracting for $D\p^{-k}$, and so we have the identity
\[
D\p^{-k}(\p^{k}(\xi_0))\cdot\textrm{det}\,D\p^{k}(\xi_0)=
\begin{pmatrix}
D_k & -B_k \\
-C_k & A_k
\end{pmatrix}
\]
Then, using the quotient rule for differentiation immediately gives
\begin{equation}\label{anglederivatives}
    \| D\t^{(k)}\| =
\biggl\|\frac{\A_{k}'\B_k-\A_k\B_{k}'}{4\A^{2}_{k}+\B^{2}_{k}}\biggr\|
\text{ and } 
\| D\t^{(k)}_{k}\| =
\biggl\|\frac{\mathcal{D}_{k}'\mathcal{C}_k-\mathcal{D}_k\mathcal{C}_{k}'}
{4\mathcal{C}^{2}_{k}+\mathcal{D}^{2}_{k}}\biggr\|.
\end{equation}

\begin{claim}\label{anglederivatives1}
      \(
      |\mathcal A_{k}|, |\mathcal B_{k}|, |\mathcal C_{k}|,
    |\mathcal D_{k}|  \leq 4\|D\varphi^{k}\|^{2}
\)
and 
\(
    \|\mathcal A'_{k}\|, \|\mathcal B'_{k}\|, |\mathcal C'_{k}\|,
    \|\mathcal D'_{k}\|  \leq 16\|D\varphi^{k}\| \|D^{2}\varphi^{k}\|
    \)
 \end{claim}
\begin{proof}
    For the first set of estimates observe that
    each partial derivative \( A_{k}, B_{k}, C_{k}, D_{k} \) of
    \( D\varphi^{k} \) is \( \leq \|D\varphi^{k}\| \). Then
    \( |\mathcal A_{k}| =  |A_k B_k+C_k D_k| \leq
    2\|D\varphi^{k}\|^{2}\). The same reasoning gives the estimates
    in the other cases.
    To estimate the derivatives, write
    \(  \|\mathcal A'_{k}\| = |A'_{k}B_{k}+A_{k}B'_{k}+
    C'_{k}D_{k}+C_{k}D'_{k}| \).  Now \( |A'_{k}| \leq 2
    \|D^{2}\varphi^{k}\| \) and similarly for the other terms.
    \end{proof}

\begin{claim}\label{anglederivatives2}
    \(
    4\mathcal{C}^{2}_{k}+\mathcal{D}^{2}_{k}=4\A^{2}_{k}+\B^{2}_{k}=
    (E^{2}_{k}-F^{2}_{k})^2.
    \)
  \end{claim}
  \begin{proof}
      Notice first of all that \( E_{k}^{2}, F_{k}^{2} \) are
      eigenvalues of
      \[
      (D\p^{k})^{T}D\p^{k} =
      \begin{pmatrix}
      A_{k} & B_{k} \\ C_{k} & D_{k}
\end{pmatrix}
 \begin{pmatrix}
      A_{k} & C_{k} \\ B_{k} & D_{k}
\end{pmatrix}
= 
 \begin{pmatrix}
      A^{2}_{k}+B_{k}^{2} & A_{k}C_{k}+ D_{k}B_{k}
      \\ A_{k}C_{k}+ D_{k}B_{k} & C^{2}_{k} + D^{2}_{k}
\end{pmatrix}
\]
In particular \( E_{k}^{2}, F_{k}^{2} \) are the two roots of the
characteristic equation
\( 
\lambda^{2}-\lambda (A_{k}^{2}+B_{k}^{2}+C_{k}^{2}+D_{k}^{2}) +
 (A^{2}_{k}+B_{k}^{2}) (C^{2}_{k} + D^{2}_{k}) -
(A_{k}C_{k} + B_{k}D_{k})^{2} = 0
\)
and therefore, by the general formula for quadratic equations, we have
\( 
F_{k}^{2}+E_{k}^{2} = A_{k}^{2}+B_{k}^{2}+C_{k}^{2}+D_{k}^{2}
\) and
\(
E_{k}^{2} F_{k}^{2} =  (A^{2}_{k}+B_{k}^{2}) (C^{2}_{k} +
D^{2}_{k}) - (A_{k}C_{k} + B_{k}D_{k})^{2}.
\)
From this one can easily check that
\(4\mathcal{C}^{2}_{k}+\mathcal{D}^{2}_{k}=4\A^{2}_{k}+\B^{2}_{k}=
    (E^{2}_{k}-F^{2}_{k})^2=(E^{2}_{k}+F^{2}_{k})^2-4E^{2}_{k}F^{2}_{k}.\)
\end{proof}

Substituting the estimates of Claims
\ref{anglederivatives1}-\ref{anglederivatives2} into
\eqref{anglederivatives} and using hyperbolicity and distortion
conditions this gives
\[ 
\|D\theta^{(k)}\|, \|D\theta^{(k)}_{k}\| \leq 128
\frac{\|D\varphi^{k}\|^{3} \|D^{2}\varphi^{k}\|}
{(E_{k}^{2}-F_{k}^{2})^{2}}
\leq 128\frac{\|D\varphi^{k}\|^{3} \|D^{2}\varphi^{k}\|}
{F^{4}_k(1-P_kQ_k H_{k}^{2})^{2}}
\leq\frac{2048}{9}\frac{\delta_{k}}{H_{k}}.
\]
To estimate
$D\t^{(k)}_{k+1}$ we write
$e^{(k)}_{k+1}=\tilde{E}_{k+1}
(\cos\t^{(k)}_{k+1},\sin\t^{(k)}_{k+1})$, so that
\begin{equation*}
\tan\t^{(k)}_{k+1}= \frac{C_1(z_k)\cos\t^{(k)}_{k}+D_1(z_k)\sin\t^{(k)}_{k}}
{A_1(z_k)\cos\t^{(k)}_{k}+B_1(z_k)\sin\t^{(k)}_{k}}=\frac{\M_k}{\N_k}.
\end{equation*}
Then
\begin{equation*}
\|D\t^{(k)}_{k+1}\|= \biggl\|\frac{\N_k\M_{k}'-\M_k\N_{k}'}
{\M^{2}_{k}+\N^{2}_{k}}\biggr\|
\quad\textrm{with}\quad
\M^{2}_{k}+\N^{2}_{k}=\frac{\|e^{(k)}_{k+1}\|^2}{\|e^{(k)}_{k}\|^2}\geq
\frac{1}{\|D\p^{-1}(z_k)\|^2}.
\end{equation*}
By inspecting this expression for $\|D\t^{(k)}_{k+1}\|$
the following bound is obtained:
\begin{equation*}
\begin{split}
\|D\t^{(k)}_{k+1}\| &\leq
2\|D\varphi^{-1}(z_k)\|^2\bigl\{2\|D\varphi(z_k)\|\bigl(
2\|D\varphi^2(z_k)\|\cdot \|D\varphi^k(z)\|+2\|D\varphi(z_k)\|
\cdot\|D\t^{(k)}_{k}\|\bigr)\bigr\}\\
&\leq 8Q^{2}_k\bigl(P^{2}_k\|D\t^{(k)}_k\|+P_k
\tilde{P}_k F_k\bigr).\\ 
\end{split} 
\end{equation*}
Putting together the estimates for $\|D\t^{(k)}_{k+1}\|$ and 
$\|D\t^{(k+1)}_{k+1}\|$, we obtain
\[
\|D\phi^{(k)}_{k+1}\|\leq \frac{2048\delta_{k+1}}{9H_{k+1}}+
8(P_kQ_k)^2\frac{\delta_{k}}{H_{k}}  
+8Q^{2}_kP_k\tilde{P}_k F_{k}. 
\] 
\end{proof}
\begin{sublemma} 
    \[ 
    \|DE_{k}\|, \|DF_{k}\| \leq \frac{2057}{9} \delta_{k} F_k/H_{k}
       \quad\text{ and } \quad
       \|D H_{k}\| \leq \frac{2066}{9} \delta_k
    \]
 \end{sublemma} 
    \begin{proof} 
We first estimate 
$D_z E_{k}=D\|e^{(k)}_{k}\|$. The corresponding estimate for $D_z F_k$ is
identical. By direct differentiation we have,
\(
D_{z} e^{(k)}_{k} =\,D^2\varphi^k(z)e^{(k)}+D\varphi^k\cdot De^{(k)}
\)
and hence by Lemma \ref{second derivative} and the estimate for $\|D\theta^{(k)}\|$
we have:
\[
\|D_{z} e^{(k)}_{k}\| \leq\,\|D^2\varphi^k(z)\|+\|D\varphi^k(z)\|\cdot
\|D_{z} e^{(k)}\|\leq\frac{\delta_kF_k}{H_k}+\frac{2048F_k\delta_k}{9H_k}
=\frac{2057\delta_kF_k}{9H_k}.
\]
Since
$D\|e^{(k)}_{k}\|=(e^{(k)}_{k}\cdot De^{(k)}_{k})\|e^{(k)}_{k}\|^{-1}$
it follows that
$\|D_z E_k\|\leq\,\|D e^{(k)}_{k}\|$ and therefore
\(  \|DE_{k}\|\leq  2057\delta_{k} F_{k}/9H_{k}\).
 Using the fact that \( \det D\varphi^{k} =
E_{k}F_{k} \) and the quotient rule for differentiation, we get
    \[ 
    DH_{k}=D\left(\frac{E_{k}}{F_{k}}\right) =
    D\left(\frac{\det D\varphi^{k}}{F_{k}^{2}}\right) =
    \frac{D(\det D\varphi^{k})}{F_{k}^{2}} - \frac{2 E_{k}
    DF_{k}}{F_{k}^{2}}.
    \]
By the estimates for $DF_k$ and Lemma \ref{second derivative} we then get
\(  \|DH_{k}\| \leq 2066\delta_{k}/9\).
    \end{proof}
To complete the proof of Lemma \ref{derivativeconvergence}, 
equations (\ref{hyp0}) and (\ref{hyp3}) give (for $k\geq k_0$):
\begin{equation*} 
     \begin{split} 
     \|D\phi^{(k)}\| &\leq 
     |H_{k+1}|(1+\tan^{2}\phi^{(k)}_{k+1})\|D\phi_{k+1}^{(k)}\|  
        + \|D H_{k+1}\| \cdot |\tan \phi^{(k)}_{k+1}|\\ 
&\leq |H_{k+1}|(1+\tan^{2}\phi^{(k)}_{k+1})\bigl(\|D\theta^{(k+1)}_{k+1}\| 
    + \|D\theta^{(k)}_{k+1}\|\bigr) \\ 
    &\qquad    + \|D H_{k+1}\| \cdot |\tan \phi^{(k)}_{k+1}|\\ 
&\leq |H_{k+1}|\bigl(1+ 4{P^{2}_kQ^{2}_k} \bigr) 
\biggl(   
\frac{2048\delta_{k+1}}{9H_{k+1}}+
8(P_kQ_k)^2\frac{\delta_{k}}{H_{k}}  
+8Q^{2}_kP_k\tilde{P}_k F_{k}
\biggr)\\ 
&\qquad+\frac{4132}{9}P_kQ_k\delta_{k+1}.\\ 
\end{split} 
\end{equation*}
 Collecting all the terms and using  \eqref{minidistortion} together with
the hyperbolicity assumptions we obtain 
\begin{equation*}
\|D\phi^{(k)}\|\leq 1597(p_kq_k)^2\delta_{k+1}+40(p_kq_k)^5\delta_{k}
+40(p_kq_k)^3q_{k}^2\tilde{p}_k\gamma^{*}_{k+1}.
\end{equation*}
This gives us the required
estimate for $\|D\phi^{(k)}\|$. To get the estimate for
$\|De^{(k)}\|$ we use the fact that 
$\|De^{(k)}\|\approx \|D\t^{(k)}\|$ with
$\t^{(k)}=\sum_{j=1}^{k-1}(\t^{(j+1)}-\t^{(j)})+\t^{(1)}.$
\end{proof}

\subsection{Global convergence}\label{sec-regularity1} 

We have seen above that the contractive directions converge pointwise
under some very mild hyperbolicity conditions.  We now want to show that
the curves \( \mathcal E^{(k)}(z) \) converge to some limit curve \(
\mathcal E^{\infty}(z) \).  
Let $z^{(k)}_t$ and $z^{(k+1)}_t$ be parametrizations by arclength of
the two curves \( \mathcal E^{(k)}(z) \) and \( \mathcal E^{(k+1)}(z)
\) with \( z^{(k)}_{0}= z^{(k+1)}_{0} = z \). 

\begin{lemma} 
     \label{leafcontraction}
For 
every \( k\geq k_{0} \) and \( t \) such that 
$z^{(k)}_t$ and $z^{(k+1)}_t$  are both defined, 
we have 
    \[  
    |z^{(k)}_{t}- z^{(k+1)}_{t}| \leq 
    t\xi_k e^{Lt}. 
    \]
\end{lemma} 

\begin{proof} 
By standard calculus  we have
      \[  
      z^{(k)}_{t} = z_{0} + \int_{0}^{t}e^{(k)}(z_{s}) ds
      \quad\text{and}\quad
         z^{(k+1)}_t = \tilde
    z_{0} + \int_{0}^{t}e^{(k+1)}( z^{(k+1)}_{s}) ds
      \]
and therefore 
      \begin{equation}\label{integral}
    |z^{(k)}_{t} - z^{(k+1)}_{t}| = \int_{0}^{t}
\|e^{(k)}(z^{(k)}_{s}) -e^{(k+1)}(z^{(k+1)}_{s}) \| ds
      \end{equation}
 By the Mean Value Theorem and Lemma \ref{angleconvergence}
 we have 
\begin{equation} 
    \label{s1} 
    \|e^{(k)}(z^{(k)}_{s}) -e^{(k)}(z^{(k+1)}_{s})\| 
    \leq 
    \|De^{(k)}\| |z^{(k)}_{s}-z^{(k+1)}_{s}| 
    \leq L|z^{(k)}_{s}-z^{(k+1)}_{s}|.
    \end{equation}
By Lemma \ref{derivativeconvergence} we have
     \begin{equation}\label{s2}
\|e^{(k)}(z^{(k+1)}_{s})
-e^{(k+1)}(z^{(k+1)}_{s})\| \leq |\phi^{(k)}| \leq \xi_k.
\end{equation}
By the triangle inequality, \eqref{s1}-\eqref{s2} give
\begin{equation}\label{uniformconvergence}
\|e^{(k)}(z^{(k)}_{s}) -e^{(k+1)}(z^{(k+1)}_{s}) \|
\leq 
L |z^{(k)}_{s}-z^{(k+1)}_{s}| + \xi_k
\end{equation}
Substituting \eqref{uniformconvergence} into
\eqref{integral}
and using Gronwall's inequality gives
 \begin{equation}\label{uniformpointwiseconvergence}
|z^{(k)}_{t}- z^{(k+1)}_{t}| \leq  t\xi_k + \int_{0}^{t}
    L|z^{(k)}_{s}-z^{(k+1)}_{s}| ds \leq
    t\xi_k e^{L t}.
\end{equation}
\end{proof}

\section{The infinite time local stable manifold}

In this section we apply the convergence estimates obtained above to
show that the local stable manifold converge to a smooth curve of
positive length and on which we have some controlled contraction
estimates. 

\subsection{Geometry}\label{geom_section}
Lemma \ref{leafcontraction} gives a bound on the distance between 
finite time local stable manifolds of different order.  However we 
have so far no guarantee that these manifolds all have some uniformly 
positive length. This depends on some delicate relationship between 
the geometry of the images of the neighbourhoods \( \mathcal N^{(k)} \) 
and the position of the finite time local stable manifolds in \( \mathcal N^{(k)}  \). 
Here we show that we can find some uniform
lower bound for the length of all local stable manifolds. 
For \( \varepsilon>0 \) and \( k \geq 1\), let 
\label{strategygeom}
\[
\w_{k}= \w_{k}(\varepsilon) = \varepsilon e^{L \varepsilon}  \xi_{k}
\]
and 
\[ 
\mathcal E^{(k)}(z, \varepsilon) = \{\xi\in \mathcal E^{(k)}(z) :
d_{\mathcal E}(\xi, z) \leq \varepsilon \}. 
\]
Recall that the constant \( L \) is determined in Lemma
\ref{derivativeconvergence}. 
Here the distance \( d_{\mathcal E}  \) is defined to be the distance 
measure inside \( \mathcal E^{(k)}(z) \) so that \( \mathcal
E^{(k)}(z,\varepsilon) \) is just a subset of \( \mathcal E^{(k)}(z) \)
which extends by a length of at most \( \varepsilon \) on both sides
of \( z \). If \( \mathcal E^{(k)}(z) \) extends by a length of less than \(
\varepsilon \) on one or both sides of \( z \) then \( \mathcal
E^{(k)}(z,\varepsilon) \) coincides with \( \mathcal E^{(k)} \) on
the corresponding sides. For simplicity we shall generally omit the \( 
\varepsilon \) from the notation and thus use the previous notation \( 
\mathcal E^{(k)}(z) \) to denote the local stable manifold of order \( 
k \) restricted to a curve of length at most \( \varepsilon \) on
either side of  \( z \). 
Let \[
\mathcal{T}_{\omega_{k}}(\mathcal E^{(k)}(z))= \{\xi: d(\xi, \mathcal 
E^{(k)}(z)) \leq \omega_{k}\}
\]
denote a 
neighbourhood of $\mathcal E^{(k)}(z)$ of size 
$\omega_k$. 
At this point we are ready to make explicit our choice of \(
\varepsilon \): we choose \( \varepsilon > 0 \)
small enough so that 
\begin{equation}\label{epsilon1}
    \varepsilon \Gamma < 1\quad\text{ and }\quad 
    e^{\varepsilon L} < 2 
\end{equation}
where \( \Gamma \) is the constant used in the definition of condition
\( (**) \), and such that 
for all  \( k\geq k_{0} \)
we have
\begin{equation}\label{epsilon}
\mathcal T_{\omega_{k}}(\mathcal E^{(k)}(z)) \subset
\mathcal N^{(k_{0})}.
\end{equation}
Notice that \eqref{epsilon} is possible because \( k_{0} \) and  \(
\mathcal N^{(k_{0})} \) are fixed and \( |\mathcal E^{(k)}(z)| \) and 
\( \omega_{k} \) can be made
arbitrarily small for \( k\geq k_{0} \) 
by taking \( \varepsilon \) small and using the fact
that \( \xi_{k}\to 0 \) by the summability condition \( (*) \). 
With this choice of \( \varepsilon \) we can then state and prove the main
result of this section. 
\begin{lemma}\label{l:epsilon}
    For all \( k\geq k_{0} \) we have 
    \[ 
    |\mathcal E^{(k)}(z)| = \varepsilon.
    \]
 \end{lemma}   
It follows that each finite time local stable manifold \( \mathcal
E^{(k)}(z) \) can be parametrized by arclength as \( z^{(k)}_{t} \)
with \( t\in [-\varepsilon, \varepsilon] \) and \( z^{(k)}_{0}=z \). 
By Lemma \ref{leafcontraction}, the pointwise limit 
\[ 
z^{(\infty)}_{t} = \lim_{k\to\infty} z_{t}^{(k)}
\]
exists for each \( t\in [-\varepsilon, \varepsilon] \) 
and defines the set 
 \[ \mathcal E^{\infty}(z) = \{z_{t}^{(\infty)}: 
 t\in [-\varepsilon, \varepsilon] \}.
     \]  
 In the following sections we will show that 
\( \mathcal E^{(\infty)}(z) \) is a smooth curve,  that 
 \(  |\mathcal E^{\infty}(z)| \geq \varepsilon \), and that it belongs to
 the stable manifolds of \( z \). 
 
 First of all we prove 
\begin{lemma}\label{l:epsilonk}
    For all \( k\geq k_{0} \) we have 
    \begin{equation}\label{epsilonk}
    \mathcal T_{\omega_{k}}(\mathcal E^{(k)}(z)) \subset
    \mathcal N^{(k + 1)}.
    \end{equation}
\end{lemma}
The proof of Lemma \ref{l:epsilonk} is a crucial step in the overall
argument and the only place in which condition \( (**) \) is used. 
Compare \eqref{epsilon} and \eqref{epsilonk}: condition
\eqref{epsilon} follows immediately by taking \( \varepsilon  \)
small, without any additional geometrical considerations. On the
other hand, \eqref{epsilonk} requires a non-trivial control of the
geometry of \( \mathcal E^{(k)}(z) \) in \( \mathcal N^{(k+1)} \). 

\begin{proof}
    We prove \eqref{epsilonk} inductively by showing that for all \(
    k\geq  j \geq k_{0} \) we have the implication
    \[ 
    \mathcal T_{\omega_{k}}(\mathcal E^{(k)}(z)) \subset
     \mathcal N^{(j)} \ \Rightarrow \ 
     \mathcal T_{\omega_{k}}(\mathcal E^{(k)}(z)) \subset
	  \mathcal N^{(j+1)}.
    \]
    Together with \eqref{epsilon}, which provides the first step of the
    induction for \( j=k_{0} \), this gives \eqref{epsilonk}. 
    Thus we need to prove that for all 
    \( x\in \mathcal T_{\omega_{k}}(\mathcal E^{(k)}(z)) \subset
     \mathcal N^{(j)}  \) we have 
    \(
    d(\varphi^{j}(x), \varphi^{j}(z)) \leq \varepsilon_{j}.
    \)
    We fix some point \( y\in \mathcal E^{(k)}(z) \) with 
    \( d( x,y) \leq \omega_{k} \). 
and write \(  d(\varphi^{j}(z), \varphi^{j}(x)) 
      \leq d(\varphi^{j}(z), \varphi^{j}(y)) + d(\varphi^{j}(y), 
      \varphi^{j}(x)) \). To estimate  
      \( d(\varphi^{j}(z), \varphi^{j}(y)) \) we use the fact that 
      both \( y \) and \( z \) are on \( \mathcal E^{(k)}(z) \) and
      that \(  \mathcal E^{(k)}(z)\) is contracting under \(
      \varphi^{j} \):
      by  \eqref{contractionclaim} on page 
      \pageref{contractionclaim}  we have (recall that \( d(z,y) \leq 
      \varepsilon \))
     \[ 
     d(\varphi^{j}(z), \varphi^{j}(y)) \leq \max_{\xi\in\mathcal
     E^{(k)}} \{\|e^{(k)}_{j}(\xi)\|\} \   d(z,y) \leq \tilde\gamma_{j}
     \varepsilon.
     \]
     To estimate  \( d(\varphi^{j}(y), \varphi^{j}(x) \) we simply use
     the fact that \( d(y,x) \leq \omega_{k} \) by assumption and in
     particular \( x \) and the line segment joining \( x \) and \( y \)
     lies entirely in \( \mathcal T_{\omega_{k}}(\mathcal E^{(k)}(z)) \)
     and therefore in \( \mathcal N^{(j)} \) by our inductive
     assumption. A relatively coarse estimate using the maximum
     possible expansion in \( \mathcal N^{(j)} \) thus gives
     \[
     d(\varphi^{j}(y), \varphi^{j}(x)) \leq \max_{\xi\in\mathcal
     N^{(j)}} \{\| F_{j}(\xi) \|\} \ d(y,x) \leq \max_{\xi\in\mathcal
     N^{(j)}} \{\|F_{j}\|\} \  \omega_{k}.
     \]
     For  \( k\geq k_{0} \) we have 
     \( \omega_{k} = \varepsilon e^{L\varepsilon} \xi_{k} 
     \leq \varepsilon e^{L\varepsilon} 2 p_{k}q_{k}\gamma_{k+1}
    \) and therefore
    \begin{equation*}
    d(\varphi^{j}(z),\varphi^{j}(x)) 
  \leq\varepsilon\tilde\gamma_j+ 2\varepsilon e^{\varepsilon L}p_{k}q_{k}
  \gamma_{k+1}
    \max_{\mathcal N^{(j)}} \{\|F_{j}\|\} = 
    \varepsilon \left(   \tilde\gamma_j+ 2 e^{\varepsilon L}p_{k}q_{k}
  \gamma_{k+1}
    \max_{\mathcal N^{(j)}} \{\|F_{j}\|\}   \right). 
    \end{equation*}  
By our choice of \( \varepsilon \) this is \( \leq \varepsilon_{j} \).
\end{proof}

  \begin{proof}[Proof of Lemma \ref{l:epsilon}] 
       By Lemma \ref{l:epsilonk},  \( e^{(k+1)} \) is defined in 
       \( \mathcal T_{\omega_{k}}(\mathcal E^{(k)}(z))  \), and therefore, 
       so is the integral leaf \( \mathcal E^{(k+1)} \). Let \( 
       z^{(k+1)}_{t} \) denote a parametrization of \( \mathcal E^{(k+1)} \) 
       by arclength with \( -t_{0}\leq t \leq t_{0} \) where \( t_{0} \) is 
       chosen maximal so that \( \{z^{(k+1)}_{t}\}_{t=-t_{0}}^{t_{0}} \subset 
       \mathcal T_{\omega_{k}}(\mathcal E^{(k)}(z))  \). We claim that 
       \( t_{0}\geq \varepsilon \), which proves the statement in the
       Lemma. Indeed, by Lemma \ref{leafcontraction} and the
       definition of \( \omega_{k} \)   we have 
   \[
   |z^{(k)}_{t}-z^{(k+1)}_{t}| \leq t\xi_k e^{L t} \leq \omega_{k}
   \] for all \( |t|\leq \varepsilon\). 
\end{proof}

\subsection{Smoothness}
\label{sec-regularity2}

We now want to study the regularity properties of
\( \mathcal E^{\infty}(z) \).

\begin{lemma}
    The curve \( \mathcal E^{(\infty)}(z) \) is \( C^{1+Lip} \) 
    with 
    \[ 
    |\mathcal E^{(\infty)}(z)| \geq \varepsilon.
    \]
The Lipschitz constant of the derivative is bounded above by
\( L \).
\end{lemma}

\begin{proof}
    From \eqref{uniformconvergence} we have
    \begin{equation}\label{uniformconvergence1}
\|e^{(k)}(z^{(k)}_{t})
    -e^{(k+1)}(z^{(k+1)}_{t}) \| \leq L
    |z^{(k)}_{t}-z^{(k+1)}_{t}| + \xi_k
\end{equation}
and from \eqref{uniformpointwiseconvergence} we have
\begin{equation}\label{uniformpointwiseconvergence1}
    |z^{(k)}_{t}-z^{(k+1)}_{t}| \leq t\xi_k e^{L t}.
\end{equation}
Thus, substituting \eqref{uniformpointwiseconvergence1}
into \eqref{uniformconvergence1}, we get
\begin{equation*}
    \|e^{(k)}(z^{(k)}_{t})
    -e^{(k+1)}(z^{(k+1)}_{t}) \| \leq
     Lt\xi_k e^{Lt}.
\end{equation*}
The uniform summability condition therefore implies that the sequence
of tangent directions \( e^{(k)}(z_{t}) \) is uniformly Cauchy in \( t \).
Thus by standard convergence results
they converge to the tangent directions of the limiting curve
\( \mathcal E^{(\infty)}(z) \) and this curve is \( C^{1} \).
To estimate the Lipschitz constant we 
let $x,x'\in\mathcal E^{\infty}(z)$ and write
\begin{equation*}
\|e^{\infty}(x)-e^{\infty}(x')\| \leq
\|e^{k}(x)-e^{\infty}(x)\|+\|e^{k}(x)-e^{k}(x')\|
+\|e^{\infty}(x')-e^{k}(x')\|
\end{equation*}
The middle term on the right hand side  is \( \leq L |x-x'| \) by the 
mean value theorem and Lemma \ref{derivativeconvergence}; the first
and last term are bounded by \( \sum_{j\geq k} |\phi^{(j)}| \leq
\sum_{j\geq k} \xi_{j} \). 
Since $\xi_k$ is summable, \( k \) is arbitrary, and $L$ uniform, the result follows.
\end{proof}

\subsection{Contraction}\label{sec-contraction}

Let \( z_{t}= z_{t}^{(\infty)} \) denote a parametrization by arclength of 
\( \mathcal E^{\infty}(z) \), with \( z_{0}=z \).

\begin{lemma}\label{contraction}
For any \(  t_{1}, t_{2}\in [\varepsilon,  -\varepsilon] \) and  \(
n\geq 1 \) we have
\[ 
|\varphi^{n}(z_{t_{1}})- \varphi^{n}(z_{t_{2}})| \leq
\tilde\gamma_n |z_{t_{1}}- z_{t_{2}}|.
\]

\end{lemma}

\begin{proof} 
    Write \( e^{(\infty)}_{n}=e^{(n)}_n+ (e^{(\infty)}_n- e^{(n)}_n) \). 
    Then 
\[ 
\|\varphi^{n}(z_{t_{1}}) - \varphi^{n}(z_{t_{2}})\| =
\int_{t_{1}}^{t_{2}}\|e^{(\infty)}_{n}\| dt = \int_{t_{1}}^{t_{2}}
\|D\varphi^{n} (e^{(n)}) + D\varphi^{n}(e^{(\infty)}- e^{(n)})\| dt.
\]
Clearly
\( \|D\varphi^{n} (e^{(n)})\| \leq\gamma_n \)
and, by  Lemma \ref{angleconvergence},
$$\|e^{(\infty)}(z)-e^{(n)}(z)\| \leq
\sum_{k\geq n} |\phi^{(k)}| 
\leq \sum_{k\geq n} \xi_k(z)\leq 2 \sum_{k\geq n}  p_{k}q_{k}\gamma_{k+1}.
$$
The definition of \( \tilde\gamma_{n} \) thus implies the statement in
the Lemma. 
\end{proof}

\subsection{Uniqueness}
\label{uniqueness}
Here we show that the local stable manifold we have constructed is
unique. That is, the set of points which remain in $\mathcal{N}^{(k)}(z)$
for all $k\geq 0$ must lie on the curve $\mathcal{E}^{(\infty)}(z)$.
\begin{lemma}
    The stable manifold through $z$ is unique in the sense that 
 \[
 \mathcal{E}^{(\infty)}(z) 
 =\bigcap_{k\geq k_0}\mathcal{N}^{(k)}(z)
 \]
\end{lemma}
\begin{proof}
 Suppose, by contradiction that there is some point 
 \(  x\in B_{\varepsilon_0}(z)  \) which belongs to 
 \( \bigcap_{k\geq k_0}\mathcal{N}^{(k)}(z) \) but not to 
 \( \mathcal{E}^{(\infty)}(z) \). We show that this point must eventually 
leave \( \bigcap_{k\geq k_0}\mathcal{N}^{(k)}(z) \). That is, there exists $j\geq 1$
such that $x\not\in\mathcal{N}^{(j)}(z)$. From the smoothness properties of 
the $e^{(k)}$ and $f^{(k)}$ vector fields, in particular their Lipschitz property, we may 
join $x$ to a 
point $\tilde{x}\in\mathcal{E}^{\infty}(z)$ via a curve $\gamma$ whose tangent 
direction has a strictly positive component in the $f^{(k)}$ direction. Hence we obtain:
\begin{equation*}
d(\varphi^{k}(z),\varphi^{k}(x)) \geq C\min_{x\in\mathcal{N}^{(k)}}\{F_k\} d(x,\tilde{x})
-\tilde\gamma_k d(z,\tilde{x}),
\end{equation*}
with $d(z,\tilde{x})<\varepsilon_0$. Since $F_k\to\infty$ as $k\to\infty$, there exists 
a $j\geq 1$ with $d(\varphi^{j}(z),\varphi^{j}(x))>\varepsilon_j$.
 
\end{proof}


\begin{bibsection}
    \begin{biblist}
	\bib{AbbMaj01}{article}{
	  author={Abbondandolo, Alberto},
	  author={Majer, Pietro},
	  title={On the Stable Manifold Theorem},
	  status={Preprint},
	  date={2001},
	}
	\bib{Ano67}{book}{
	  author={Anosov, D. V.},
	  title={Geodesic flows on closed Riemann manifolds with negative curvature. },
	  series={Proceedings of the Steklov Institute of Mathematics, No. 90 (1967). Translated from the Russian by S. Feder},
	  publisher={American Mathematical Society},
	  place={Providence, R.I.},
	  date={1969},
	  pages={iv+235},
	}
	\bib{BenCar91}{article}{
	  author={Benedicks, M.},
	  author={Carleson, L.},
	  title={The dynamics of the He\'non map},
	  date={1991},
	  journal={Ann. of Math.},
	  volume={133},
	  pages={73\ndash 169},
	}
	\bib{BenYou93}{article}{
	  author={Benedicks, Michael},
	  author={Young, Lai-Sang},
	  title={Sina\u \i -Bowen-Ruelle measures for certain H\'enon maps},
	  journal={Invent. Math.},
	  volume={112},
	  date={1993},
	  number={3},
	  pages={541\ndash 576},
	  issn={0020-9910},
	}
	\bib{BogMit61}{book}{
	  author={Bogoliubov, N. N.},
	  author={Mitropolsky, Y. A.},
	  title={Asymptotic methods in the theory of non-linear oscillations},
	  series={Translated from the second revised Russian edition. International Monographs on Advanced Mathematics and Physics},
	  publisher={Hindustan Publishing Corp.},
	  place={Delhi, Gordon and Breach Science Publishers, New York},
	  date={1961},
	  pages={x+537},
	}
	\bib{BogMit63}{article}{
	  author={Bogoliubov, N. N.},
	  author={Mitropolsky, Yu. A.},
	  title={The method of integral manifolds in nonlinear mechanics},
	  journal={Contributions to Differential Equations},
	  volume={2},
	  date={1963},
	  pages={123\ndash 196 (1963)},
	}
	\bib{Cha01}{article}{
	  author={Chaperon, Marc},
	  title={Some results on stable manifolds},
	  language={English, with English and French summaries},
	  journal={C. R. Acad. Sci. Paris S\'er. I Math.},
	  volume={333},
	  date={2001},
	  number={2},
	  pages={119\ndash 124},
	  issn={0764-4442},
	}
	\bib{Cha02}{article}{
	  author={Chaperon, M.},
	  title={Invariant manifolds revisited},
	  note={Dedicated to the 80th annniversary of Academician Evgeni\u \i \ Frolovich Mishchenko (Russian) (Suzdal, 2000)},
	  journal={Tr. Mat. Inst. Steklova},
	  volume={236},
	  date={2002},
	  number={Differ. Uravn. i Din. Sist.},
	  pages={428\ndash 446},
	  issn={0371-9685},
	}
	\bib{Dil60}{article}{
	  author={Diliberto, Stephen P.},
	  title={Perturbation theorems for periodic surfaces. I. Definitions and main theorems},
	  journal={Rend. Circ. Mat. Palermo (2)},
	  volume={9},
	  date={1960},
	  pages={265\ndash 299},
	}
	\bib{Dil61}{article}{
	  author={Diliberto, Stephen P.},
	  title={Perturbation theorems for periodic surfaces. II},
	  journal={Rend. Circ. Mat. Palermo (2)},
	  volume={10},
	  date={1961},
	  pages={111\ndash 161},
	}
	\bib{Fen71}{article}{
	  author={Fenichel, N},
	  title={Persistence and smoothness of invariant manifolds of flows},
	  journal={Indiana Univ.Math.J.},
	  year={1971},
	  pages={193-226},
	  volume={21},
	}
	\bib{Fen74}{article}{
	  author={Fenichel, N},
	  title={Asymptotic stability with rate conditions},
	  volume={23},
	  pages={1109-1137},
	  journal={Indiana Univ.Math.J.},
	  year={1974},
	}
	\bib{Had01s}{article}{
	  author={Hadamard, J. },
	  title={Selecta},
	  date={1901},
	  pages={163-166},
	}
	\bib{Had01}{article}{
	  author={Hadamard, J. },
	  title={Sur l'iteration et les solutions asymptotetiques des equations differentielles},
	  journal={Bull.Soc.Math. France},
	  year={1901},
	  volume={29},
	  pages={224-228},
	}
	\bib{Hal60}{article}{
	  author={Hale, Jack K.},
	  title={On the method of Krylov-Bogoliubov-Mitropolski for the existence of integral manifolds of perturbed differential systems},
	  journal={Bol. Soc. Mat. Mexicana (2)},
	  volume={5},
	  date={1960},
	  pages={51\ndash 57},
	}
	\bib{Hal61}{article}{
	  author={Hale, Jack K.},
	  title={Integral manifolds of perturbed differential systems},
	  journal={Ann. of Math. (2)},
	  volume={73},
	  date={1961},
	  pages={496\ndash 531},
	}
	\bib{HirPalPugShu70}{article}{
	  author={Hirsch, M.},
	  author={Palis, J.},
	  author={Pugh, C.},
	  author={Shub, M.},
	  title={Neighborhoods of hyperbolic sets},
	  journal={Invent. Math.},
	  volume={9},
	  date={1969/1970},
	  pages={121\ndash 134},
	}
	\bib{HirPug69}{article}{
	  author={Hirsch, Morris W.},
	  author={Pugh, Charles C.},
	  title={Stable manifolds for hyperbolic sets},
	  journal={Bull. Amer. Math. Soc.},
	  volume={75},
	  date={1969},
	  pages={149\ndash 152},
	}
	\bib{HirPug70}{article}{
	  author={Hirsch, Morris W.},
	  author={Pugh, Charles C.},
	  title={Stable manifolds and hyperbolic sets},
	  booktitle={Global Analysis (Proc. Sympos. Pure Math., Vol. XIV, Berkeley, Calif., 1968)},
	  pages={133\ndash 163},
	  publisher={Amer. Math. Soc.},
	  place={Providence, R.I.},
	  date={1970},
	}
	\bib{HirPugShu77}{article}{
	  author={Hirsch, Morris},
	  author={Pugh, Charles},
	  author={Shub, Michael},
	  title={Invariant Manifolds},
	  booktitle={Lecture Notes in Math.},
	  volume={583},
	  publisher={Springer},
	  year={1977},
	}
	\bib{HolLuz}{article}{
	  author={Holland, Mark},
	  author={Luzzatto, Stefano},
	  title={Dynamics of two dimensional maps with criticalities and singularities},
	  status={Work in progress},
	  date={2003},
	}
	\bib{Irw70}{article}{
	  author={Irwin, M.C.},
	  title={On the stable manifold theorem},
	  journal={Bull.London.Math.Soc.},
	  year={1970},
	  volume={2},
	  pages={196-198},
	}
	\bib{KatHas94}{article}{
	  author={Katok, Anatole},
	  author={Hasselblatt, Boris},
	  title={Introduction to the modern Theory of Smooth Dynamical Systems},
	  publisher={Cambridge University Press},
	  year={1994},
	  place={Cambridge, New York},
	}
	\bib{Kel67}{article}{
	  author={Kelley, Al},
	  title={The stable, center-stable, center, center-unstable, unstable manifolds},
	  journal={J. Differential Equations},
	  volume={3},
	  date={1967},
	  pages={546\ndash 570},
	}
	\bib{LuzVia2}{article}{
	  author={\href {http://www.ma.ic.ac.uk/~luzzatto}{Stefano Luzzatto}},
	  author={Viana, Marcelo},
	  title={Lorenz-like attractors without continuous invariant foliations},
	  status={Preprint},
	  year={2003},
	}
	\bib{MorVia93}{article}{
	  author={Mora, Leonardo},
	  author={Viana, Marcelo},
	  title={Abundance of strange attractors},
	  journal={Acta Math.},
	  volume={171},
	  date={1993},
	  number={1},
	  pages={1\ndash 71},
	  issn={0001-5962},
	}
	\bib{Ose68}{article}{
	  author={Oseledec, V. I.},
	  title={A multiplicative ergodic theorem. Characteristic Ljapunov, exponents of dynamical systems},
	  journal={Transactions of the Moscow Mathematical Society},
	  date={1968},
	  volume={19},
	  publisher={American Mathematical Society},
	  place={Providence, R.I.},
	}
	\bib{Per28}{article}{
	  author={Perron, Oskar},
	  title={\"Uber Stabilit\"at und asymptotisches Verhalten der Integrale von Differentialgleichungssytemen},
	  journal={Math. Z},
	  date={1929},
	  volume={161},
	  pages={41--64},
	}
	\bib{Per30}{article}{
	  author={Perron, Oskar},
	  title={Die Stabilit\"atsfrage bei Differentialgleichungen},
	  journal={Math. Z},
	  volume={32},
	  date={1930},
	  pages={703-728},
	}
	\bib{Pes76}{article}{
	  title={Families of invariant manifolds corresponding to non-zero characteristic exponents},
	  author={Ya. Pesin},
	  journal={Math. USSR. Izv.},
	  volume={10},
	  pages={1261\ndash 1302},
	  year={1976},
	}
	\bib{Pes77}{article}{
	  title={Characteristic Lyapunov exponents and smooth ergodic theory},
	  author={Pesin, Ya. B. },
	  journal={Russian Math. Surveys},
	  volume={324},
	  pages={55\ndash 114},
	  year={1977},
	}
	\bib{Poi86s}{article}{
	  author={Poincar\'e, Henri},
	  title={Oeuvres},
	  volume={1},
	  date={1886},
	  pages={202--204},
	}
	\bib{Sac67}{article}{
	  author={Sacker, Robert J.},
	  title={A perturbation theorem for invariant Riemannian manifolds},
	  booktitle={Differential Equations and Dynamical Systems (Proc. Internat. Sympos., Mayaguez, P.R., 1965)},
	  pages={43\ndash 54},
	  publisher={Academic Press},
	  place={New York},
	  date={1967},
	}
	\bib{Sma67}{article}{
	  author={Smale, S.},
	  title={Differentiable dynamical systems},
	  journal={Bull. Amer. Math. Soc.},
	  volume={73},
	  date={1967},
	  pages={747\ndash 817},
	}
	\bib{Ste55}{article}{
	  author={Sternberg, Shlomo},
	  title={On the behavior of invariant curves near a hyperbolic point of asurface transformation},
	  journal={Amer. J. Math.},
	  volume={77},
	  date={1955},
	  pages={526\ndash 534},
	}
	\bib{WanYou01}{article}{
	  author={Wang, Qiudong},
	  author={Young, Lai-Sang},
	  title={Strange attractors with one direction of instability},
	  journal={Comm. Math. Phys.},
	  volume={218},
	  date={2001},
	  number={1},
	  pages={1\ndash 97},
	}

    \end{biblist}
 \end{bibsection}
 \end{document}